\documentclass[english]{article}
\usepackage[T1]{fontenc}
\usepackage[latin9]{inputenc}
\usepackage{amssymb}
\usepackage{esint}
\usepackage{babel}
\begin{document}

\title{Generalization of Lambert W-function, Bessel polynomials and transcendental
equations}

\author{Giorgio Mugnaini %
\thanks{\texttt{email: giorgiomugnaini1974@gmail.com}%
} }

\date{26/12/2014}
\maketitle
\begin{abstract}
Employing the Lagrange inverting series, a solution of a mixed exponential/quadratic
transcendental equation, that can be considered a generalization of
the equation defining Lambert W-function, has been found in terms
of Bessel orthogonal polynomials, and a novel representation for Bessel
polynomials has been found. Some other cases of transcendental equations
formally solved by means of classical orthogonal polynomials have
been discussed, suggesting a link between Rodrigues formulas and the
terms of Lagrange series to be investigated.
\end{abstract}

\section{Lambert-like transcendental equation}

We consider the following transcendental equation: 
\begin{equation}
(x-a)(x-b)=le^{x}
\end{equation}
\label{eq1}

This transcendental equation emerges in the study of the electron
in the potential of $H{}_{2}^{+}$ion (see \cite{Mezo14}) , and it
can be considered a generalization of Lambert W-function defining
equation. Indeed, in the limit case -b= l =$\infty$we obtain:

\[
x-a=le^{x}
\]

that can be solved by means of Lambert W function:

\[
x=a-W(-le^{a})
\]

Now we will produce a solution, based on Lagrange inversion series,
of equation \ref{eq1}.

We can rewrite such equation as follows:

\[
x=a+\frac{le^{x}}{x-b}
\]

where only trivial algebraic manipulations are employed.

Now, remembering that the solution of the equation :

\[
x=a+lf(x)
\]

can be written by means of the Lagrange inversion:

\[
x=a+\sum_{n=1}^{\infty}\frac{l^{n}}{n!}\left[\left(\frac{d}{dx}\right)^{n-1}f\left(x\right)^{n}\right]_{x=a}
\]

we obtain a formal series solving equation \ref{eq1}:

\[
x=a+\sum_{n=1}^{\infty}\frac{l^{n}}{n!}\left[\left(\frac{d}{dx}\right)^{n-1}\left(\frac{e^{x}}{x-b}\right)^{n}\right]_{x=a}
\]

Developing the derivative:

\[
x=a+\sum_{n=1}^{\infty}\frac{l^{n}}{n!}\left[\sum_{k=0}^{n-1}\left(\begin{array}{c}
n-1\\
k
\end{array}\right)\left[\left(\frac{d}{dx}\right)^{n-1-k}e^{nx}\right]\left[\left(\frac{d}{dx}\right)^{k}\left(\frac{1}{x-b}\right)^{n}\right]\right]_{x=a}
\]

and remembering that:

\[
\left(\frac{d}{dx}\right)^{k}\left(\frac{1}{x-b}\right)^{n}=\frac{(n+k-1)!}{(n-1)!}\left(\frac{1}{x-b}\right)^{n+k}(-1)^{k}
\]

substituting, we found:

\[
x=a+\sum_{n=1}^{\infty}\frac{l^{n}}{n!}\sum_{0}^{n-1}\left(\begin{array}{c}
n-1\\
k
\end{array}\right)n^{n-1-k}e^{na}(-1)^{k}\frac{(n+k-1)!}{(n-1)!}\left(\frac{1}{a-b}\right)^{n+k}
\]

\[
x=a+\sum_{n=1}^{\infty}\frac{l^{n}}{n!}\sum_{0}^{n-1}\frac{(n-1)!}{(n-1-k)!k!}n^{n-1-k}e^{na}(-1)^{k}\frac{(n+k-1)!}{(n-1)!}\left(\frac{1}{a-b}\right)^{n+k}
\]

grouping terms not depending by k to the left of second summation
sign:

\[
x=a+\sum_{n=1}^{\infty}\frac{1}{n!n}\left(\frac{nle^{a}}{a-b}\right)^{n}\sum_{0}^{n-1}\frac{(n-1+k)!}{(n-1-k)!k!}\left(\frac{-1}{n(a-b)}\right)^{k}
\]

Remembering the definition of Bessel polynomials (see \cite{Krall49})
:

\[
B_{n}(z)=\sum_{k=0}^{n}\frac{(n+k)!}{(n-k)!k!}\left(\frac{z}{2}\right)^{k}
\]

we can rewrite the solution in a more compact form:

\begin{equation}
x=a+\sum_{n=1}^{\infty}\frac{1}{n!n}\left(\frac{nle^{a}}{a-b}\right)^{n}B_{n-1}\left(\frac{-2}{n(a-b)}\right)
\end{equation}

Symmetrically we can obtain another solution series by swap of parameters
a and b:

\begin{equation}
x=b+\sum_{n=1}^{\infty}\frac{1}{n!n}\left(-\frac{nle^{b}}{a-b}\right)^{n}B_{n-1}\left(\frac{+2}{n(a-b)}\right)
\end{equation}

\section{A novel representation for Bessel polynomials}

It is interesting to observe that, in the development of above seen
calculations, we have encountered a novel representation for Bessel
polynomials:

\begin{equation}
B_{n-1}\left(\frac{-2}{x}\right)=e^{-x}x^{n}\left(\frac{d}{dx}\right)^{n-1}\left[\frac{e^{x}}{x^{n}}\right]=\sum_{k=0}^{n}\frac{(n-1+k)!}{(n-1-k)!(k)!}\left(-\frac{1}{x}\right)^{k}\label{eq:novel}
\end{equation}

The first cases for n=0..5, of the polynomials of type$e^{-x}x^{n+1}\left(\frac{d}{dx}\right)^{n}\left[\frac{e^{x}}{x^{n+1}}\right]$are:

\[
e^{-x}x\left(\frac{d}{dx}\right)^{0}\left[\frac{e^{x}}{x^{1}}\right]=1
\]

\[
e^{-x}x^{2}\left(\frac{d}{dx}\right)^{1}\left[\frac{e^{x}}{x^{2}}\right]=1-\frac{2}{x}
\]

\[
e^{-x}x^{3}\left(\frac{d}{dx}\right)^{2}\left[\frac{e^{x}}{x^{3}}\right]=1-\frac{6}{x}+\frac{12}{x^{2}}
\]

\[
e^{-x}x^{4}\left(\frac{d}{dx}\right)^{3}\left[\frac{e^{x}}{x^{4}}\right]=1-\frac{12}{x}+\frac{60}{x^{2}}-\frac{120}{x^{3}}
\]

\[
e^{-x}x^{5}\left(\frac{d}{dx}\right)^{4}\left[\frac{e^{x}}{x^{5}}\right]=1-\frac{20}{x}+\frac{180}{x^{2}}-\frac{840}{x^{3}}+\frac{1680}{x^{4}}
\]

\[
e^{-x}x^{6}\left(\frac{d}{dx}\right)^{5}\left[\frac{e^{x}}{x^{6}}\right]=1-\frac{30}{x}+\frac{420}{x^{2}}-\frac{3360}{x^{3}}+\frac{15120}{x^{4}}-\frac{30240}{x^{5}}
\]

whereas the classical Rodrigues representation for Bessel polynomials
gives :

\[
e^{-\frac{1}{x}}x\left(\frac{d}{dx}\right)^{0}\left[e^{\frac{1}{z}}x^{2*0}\right]=(-1)^{0}1
\]

\[
e^{-\frac{1}{x}}\left(\frac{d}{dx}\right)^{1}\left[e^{\frac{1}{z}}x^{2*1}\right]=(-1)^{1}\left(1-2x\right)
\]

\[
e^{-\frac{1}{x}}\left(\frac{d}{dx}\right)^{2}\left[e^{\frac{1}{z}}x^{2*2}\right]=(-1)^{2}\left(1-6x+12x^{2}\right)
\]

\[
e^{-\frac{1}{x}}\left(\frac{d}{dx}\right)^{3}\left[e^{\frac{1}{z}}x^{3*2}\right]=(-1)^{3}\left(1-12x+60x^{2}-120x^{3}\right)
\]

\[
e^{-\frac{1}{x}}\left(\frac{d}{dx}\right)^{4}\left[e^{\frac{1}{z}}x^{4*2}\right]=(-1)^{4}\left(1-20x+180x^{2}-840x^{3}+1680x^{4}\right)
\]

\[
e^{-\frac{1}{x}}\left(\frac{d}{dx}\right)^{5}\left[e^{\frac{1}{z}}x^{5*2}\right]=(-1)^{5}\left(1-30x+420x^{2}-3360x^{3}+15120x^{4}-30240x^{5}\right)
\]

A direct derivation of the proposed representation of Bessel polynomials
in $\frac{1}{x}$, obtained from the classical Rodrigues formula (\cite{Krall49})
will be provided in the next section.

\section{Rodrigues formula and transcendental equations}

The appearance of orthogonal polynomials (i.e. Bessel polynomials)
in solution of above seen transcendental equation is not surprising,
indeed remembering the generalized Rodrigues formula (\cite{Hochstadt}):

\[
P_{n}(x)=\frac{1}{W(x)}\left(\frac{d}{dx}\right)^{n}\left[W(x)Q(x)^{n}\right]
\]

it can be placed in a form resembling the terms of the Lagrange inversion
series, by suitable choose of terms Q and W:

\[
P_{n}(x)=\left(\frac{d}{dx}\right)^{n-1}f(x)^{n}
\]

when f(x) is a product of suitable rational function and exponential
in variable x,

In following subsections, we will examine some examples of interpretation
of Lagrange series terms as Rodrigues formulas.

\subsection{Equation $\frac{x-s}{x-t}=le^{x}$ and Laguerre polynomials}

We can remember\cite{Mezo14} for an employment of Laguerre polynomials
for solution of a similar transcendental equation (another generalization
of original Lambert W equation):

\begin{equation}
\frac{x-s}{x-t}=le^{x}
\end{equation}

which can be formally solved by:

\begin{equation}
x=t+\sum_{n=1}^{\infty}\frac{(t-s)l^{n}}{n}L{}_{n-1}^{(1)}\left(n(t-s)\right)\label{eq:laguerre_series}
\end{equation}

where L (x) denotes the generalized Laguerre polynomials (in\cite{Mezo14},
an equivalent formula based on first derivative of Laguerre polynomials
was found).

In the mentioned paper an non-trivial induction proof was given, but
a more straightforward derivation can be obtained if we recognize
that the above equation can be rewritten as:

\begin{equation}
x-s=le^{x}\left(x-t\right)
\end{equation}

that can be solved by Lagrange inversion, obtaining for n-th term
in the series (having placed: $f(x)=e^{x}\left(x-t\right)$):

\begin{equation}
P_{n}(x)=\left(\frac{d}{dx}\right)^{n-1}e^{nx}\left(x-t\right)^{n}
\end{equation}

remembering the definition of Laguerre polynomials, the equation (\ref{eq:laguerre_series})
follows.

\subsection{Equation $x=a+le^{-x^{2}/2}$and Hermite polynomials}

The equation:

\begin{equation}
x=a+le^{-x^{2}/2}\label{eq:equazione_hermite}
\end{equation}

leads to Lagrange term :

\begin{equation}
P_{n}(x)=\left(\frac{d}{dx}\right)^{n-1}e^{-x^{2}/2}\label{eq:risolvente_hermite}
\end{equation}

remembering Rodrigues formula for Hermite polynomial (see \cite{Hochstadt}):

\begin{equation}
H_{n}(x)=e^{x^{2}/2}\left(\frac{d}{dx}\right)^{n}e^{-x^{2}/2}
\end{equation}

we can rewrite the Lagrange series term (\ref{eq:risolvente_hermite})
by means of Rodrigues formula as:

\begin{equation}
P_{n}(x)=e^{nx^{2}/2}H_{n-1}\left(\sqrt{n}x\right)
\end{equation}

therefore the complete solution of equation (\ref{eq:equazione_hermite}):

\begin{equation}
x=a+\sum_{n=1}^{\infty}\frac{l^{n}}{n!}e^{na^{2}/2}H_{n-1}\left(\sqrt{n}a\right)\label{eq:laguerre_series-1}
\end{equation}

\subsection{Equation $x=a+le^{e^{x}}$and Touchard polynomials}

The equation:

\begin{equation}
x=a+le^{e^{x}}\label{eq:equazione_touchard}
\end{equation}

leads to Lagrange term:

\begin{equation}
P_{n}(x)=\left(\frac{d}{dx}\right)^{n-1}e^{e^{x}}\label{eq:risolvente_touchard}
\end{equation}

remembering Rodrigues formula for Touchard polynomials:

\[
T_{n}(x)=e^{-e^{x}}\left(\frac{d}{dx}\right)^{n}e^{e^{x}}
\]

therefore:
\[
P_{n}(x)=e^{e^{x}}T_{n-1}\left(x+log(n)\right)
\]

and the formal solution of (\ref{eq:equazione_touchard}):

\begin{equation}
x=a+\sum_{n=1}^{\infty}\frac{l^{n}}{n!}e^{e^{a}}T_{n-1}\left(a+log(n)\right)\label{eq:laguerre_series-1-1}
\end{equation}

\subsection{Equation $x=a+lx^{2}e^{-\frac{2}{x}}$and Bessel polynomials}

The Bessel polynomials can be generated by the following Rodrigues
formula (see\cite{Krall49}):

\begin{equation}
B_{n}(x)=2^{-n}e^{\frac{2}{x}}\left(\frac{d}{dx}\right)^{n}\left[e^{-\frac{2}{x}}x^{2n}\right]
\end{equation}

An employment the Rodrigues formula for Bessel polynomials for solving%
\footnote{Replacing $x\rightarrow\frac{1}{x}$this equation becomes $x=ax^{2}+le^{-2x}$,
having same form of \ref{eq1}%
} seems natural:

\begin{equation}
x=a+lx^{2}e^{-\frac{2}{x}}\label{eq:equazione_reciproca}
\end{equation}

therefore the formal n-th term in the Lagrange solution series is:

\begin{equation}
\left(\frac{d}{dx}\right)^{n-1}\left[e^{-\frac{2}{x}}x^{2n}\right]
\end{equation}

that can be formally written in terms of Bessel polynomials as follows:$\int B_{n}(x)2^{n}e^{-\frac{2}{x}}dx$.
Therefore the complete solution of equation \ref{eq:equazione_reciproca}
becomes:

\begin{equation}
x=a+\sum_{n=1}\frac{l^{n}}{n!}\int B_{n}(a)2^{n}e^{-\frac{2}{a}}dx
\end{equation}

For the case in study equation (\ref{eq1}), we have to place:

\[
f(x)=\frac{e^{x}}{x-b}
\]

\subsection{Equation $(x-a)(x-b)=le^{x}$and Bessel polynomials in novel Rodrigues-like
form}

In previous section, starting from Lagrange series, we have found
a novel Rodrigues-like representation of Bessel polynomials:

\[
B_{n-1}\left(\frac{-2}{x}\right)=\left(-2\right)^{n}x^{n}e^{-x}\left(\frac{d}{dx}\right)^{n-1}\left[\frac{e^{x}}{x^{n}}\right]
\]

We now will provide an alternative proof; starting from the equation
(6) we change x->2/x:

\begin{equation}
B_{n}\left(\frac{-2}{x}\right)=2^{-n}e^{-x}\left(\frac{d}{d\frac{-2}{x}}\right)^{n}\left[e^{x}\frac{2^{2n}}{x^{2n}}\right]\label{eq:reciproco_deriv}
\end{equation}

\begin{equation}
B_{n}\left(\frac{-2}{x}\right)=e^{-x}\left(x^{2}\frac{d}{dx}\right)^{n}\left[\frac{e^{x}}{x^{2n}}\right]
\end{equation}

Now we will make use of the following identity (we will provide a
proof in the next section):

\begin{equation}
\left(-\frac{d}{d\frac{1}{x}}\right)^{n}=\left(x^{2}\frac{d}{dx}\right)^{n}=x^{n+1}\left(\frac{d}{dx}\right)^{n}x^{n-1}
\end{equation}

so we can rewrite equation (\ref{eq:reciproco_deriv}):

\begin{equation}
B_{n}\left(\frac{-2}{x}\right)=e^{-x}x^{n+1}\left(\frac{d}{dx}\right)^{n}\left[\frac{e^{x}}{x^{n+1}}\right]
\end{equation}

that matches (\ref{eq:novel}), after rescaling $x\rightarrow x/n$.

\section{Proof of identity $\left(x^{2}D\right)^{n}=x^{n+1}D^{n}x^{n-1}$
by induction}

Now we will provide a proof for identity :

\begin{equation}
\left(-\frac{d}{d\frac{1}{x}}\right)^{n}=\left(x^{2}D\right)^{n}=x^{n+1}D^{n}x^{n-1}\label{eq:differential}
\end{equation}

where we have placed: $D=\frac{d}{dx}$in order to save space.

We start observing that, for n=1 the following is trivial that:

\begin{equation}
\left(x^{2}D\right)=x^{n+1}D^{n}x^{n-1},n=1
\end{equation}

Therefore for n=1 the identity holds. For induction we proof that,
if the identity holds for n, then it holds also for n+1. Applying
the operator$\left(x^{2}D\right)$ 

to both members of (\ref{eq:differential}):

\begin{equation}
\left(x^{2}D\right)\left(x^{2}D\right)^{n}=\left(x^{2}D\right)\left(x^{n+1}D^{n}x^{n-1}\right)
\end{equation}

\begin{equation}
\left(x^{2}D\right)^{n+1}=x^{2}Dxx^{n}D^{n}x^{n-1}=x^{2}Dx^{n}xD^{n}x^{n-1}
\end{equation}

remembering the following commutation rules:

\begin{equation}
Dx^{n}\centerdot-x^{n}D=nx^{n-1}
\end{equation}

\begin{equation}
D^{n}x\centerdot-xD^{n}=nD^{n-1}
\end{equation}

we can write:

\begin{equation}
\left(x^{2}D\right)\left(x^{2}D\right)^{n}=x^{2}\left[Dx^{n}\centerdot\right]\left[xD^{n}\right]x^{n-1}=x^{2}\left[x^{n}D+nx^{n-1}\right]\left[D^{n}x\centerdot-nD^{n-1}\right]x^{n-1}
\end{equation}
\begin{equation}
=x^{n+2}D^{n+1}x^{n}\centerdot+x^{2}\left[nx^{n-1}D^{n}x\centerdot-x^{n}DnD^{n-1}-nx^{n-1}nD^{n-1}\right]x^{n-1}
\end{equation}

Now we can show that the term in square parentheses is null, indeed
applying above seen commutation rules:

\begin{equation}
=x^{n+2}D^{n+1}x^{n}\centerdot+x^{2}\left[nx^{n-1}xD^{n}+nx^{n-1}nD^{n-1}-nx^{n}D{}^{n}-nx^{n-1}nD^{n-1}\right]x^{n-1}=
\end{equation}

\begin{equation}
=x^{n+2}D^{n+1}x^{n}\centerdot+x^{2}\left[nx^{n}D^{n}-nx^{n}D{}^{n}+nx^{n-1}nD^{n-1}-nx^{n-1}nD^{n-1}\right]x^{n-1}=x^{n+2}D^{n+1}x^{n}\centerdot
\end{equation}

Therefore:

\begin{equation}
\left(x^{2}D\right)^{n+1}=\left(x^{2}D\right)\left(x^{2}D\right)^{n}=x^{n+2}D^{n+1}x^{n}\centerdot
\end{equation}

\section{Concluding remarks}

Employing the Lagrange inverting series, a solution of the transcendental
equation $(x-a)(x-b)=le^{x}$(that can be considered a quadratic generalization
of the equation defining Lambert W-function) has been found in terms
of Bessel orthogonal polynomials. 

Once again (see \cite{Mezo14}), a transcendental equation can be
formally solved by means of classic orthogonal polynomials, suggesting
a link between Rodrigues formulas and the terms of Lagrange series.
Moreover a novel representation for Bessel polynomials $B_{n}\left(\frac{-2}{x}\right)=e^{-x}x^{n+1}\left(\frac{d}{dx}\right)^{n}\left[\frac{e^{x}}{x^{n+1}}\right]$is
shown. Finally, an interesting and previously unknown differential
identity $\left(-\frac{d}{d\frac{1}{x}}\right)^{n}=\left(x^{2}D\right)^{n}=x^{n+1}D^{n}x^{n-1}\centerdot$
for the reciprocal differentiation has been found. 

Further investigations on radius of convergence and numerical feasibility
of above seen series (possibly by series acceleration techniques)
are required.

Moreover the occurrence of classical orthogonal polynomials in solution
of mixed exponential/polynomial transcendental equations suggests
to investigate for more general series of hypergeometric functions
generated by Lagrange inversion of transcendental equations.

\end{document}